\documentclass[12pt]{article}

\usepackage{amsmath,amssymb,enumerate,bbm}
\newcommand{\assign}{:=}
\newcommand{\nin}{\not\in}
\newcommand{\tmmathbf}[1]{\ensuremath{\boldsymbol{#1}}}
\newcommand{\tmop}[1]{\ensuremath{\operatorname{#1}}}
\newtheorem{theorem}{Theorem}[section]
\newtheorem{corollary}[theorem]{Corollary}
\newtheorem{lemma}[theorem]{Lemma}

\numberwithin{equation}{section}
\newenvironment{proof}{\noindent\textbf{Proof\ }}{\hspace*{\fill}$\Box$\medskip}

\begin{document}

\title{Distribution of determinant of matrices with restricted entries over
finite fields}\author{Le Anh Vinh\\
Mathematics Department\\
Harvard University\\
Cambridge, MA 02138\\
vinh@math.harvard.edu}\maketitle

\begin{abstract}
  For a prime power $q$, we study the distribution of determinent of matrices with restricted entries over a finite field $\mathbbm{F}_q$ of $q$ elements. More precisely, let $N_d
(\mathcal{A}; t)$ be the number of $d \times d$ matrices with entries in
$\mathcal{A}$ having determinant $t$. We show that
\[ N_d (\mathcal{A}; t) = (1 + o (1)) \frac{|\mathcal{A}|^{d^2}}{q}, \]
  if $|\mathcal{A}| = \omega(q^{\frac{d}{2d-1}})$, $d\geqslant 4$. When $q$ is a prime and $\mathcal{A}$ is a symmetric interval $[-H,H]$, we get the same result for $d\geqslant 3$. This improves a result of Ahmadi and Shparlinski (2007).
\end{abstract}

\begin{center}
\end{center}

\section{Introduction}

Throughout the paper, let $q=p^r$ where $p$ is an odd prime and $r$ is a positive integer. Let $\mathbbm{F}_q$ be a finite field of $q$ elements. The prime base field $\mathbbm{F}_p$ of $\mathbbm{F}_q$ may then be naturally identified with $\mathbbm{Z}_p$. For integer numbers $m$ and $n$, let $\mathcal{M}_{m, n} (\mathcal{A})$ denote
the set of $m \times n$ matrices with components in the set $\mathcal{A}$. In \cite{igor}, Ahmadi and Shparlinski studied some natural classes of matrices over a finite
field $\mathbbm{F}_p$ of $p$ elements ($p$ is a prime) with components in a
given subinterval $[- H, H] \subseteq [- (p - 1) / 2, (p - 1) / 2]$. Let $N_d
(\mathcal{A}; t)$ be the number of $d \times d$ matrices with entries in
$\mathcal{A}$ having determinant $t$. Ahmadi and Shparlinski \cite{igor} proved the following result (see \cite{igor} and the references therein for the motivation and related results).

\begin{theorem}\label{t1}
  (\cite[Theorem 11]{igor}) For $1 \leqslant H \leqslant (p - 1) / 2$ and $t \in
  \mathbbm{F}_p^{*}$, we have
  \[ N_d ([- H, H] ; t) = \frac{(2 H + 1)^{d^2}}{p} + O (H^{d^2 - 2} p^{1 / 2}
     (\log p)^2) . \]
\end{theorem}

Note that the proof of Theorem 11 in \cite{igor} is given only in the case $t = 1$, but
it goes through without any essential changes for arbitrary $t \in
  \mathbbm{F}_p^{*}$. The
bound of Theorem \ref{t1} is nontrivial if $H \gg p^{3 / 4 + \epsilon}$. In the case $d = 2$,
they obtained a stronger result.

\begin{theorem} \label{t2}
  (\cite[Theorem 12]{igor}) For $1 \leqslant H \leqslant (p - 1) / 2$ and $t \in
  \mathbbm{F}_p^{*}$, we have
  \[ N_2 ([- H, H] ; t) = \frac{(2 H + 1)^4}{p} + O (H^2 p^{o (1)}) . \]
\end{theorem}

Again, the proof of Theorem 12 in \cite{igor} is given only in the case $t = 1$, but it
goes through without any changes for arbitrary $t \in \mathbbm{F}_p^{*}$. The
bound of Theorem \ref{t2} is nontrivial if $H \gg p^{1 / 2+\epsilon}$. 

Covert et al. \cite{covert} studied this problem in a more general setting. More
precisely, define vol$(\tmmathbf{x}^1, \ldots, \tmmathbf{x}^d)$ to be the
determinant of the matrix whose rows are $\tmmathbf{x}^j$s. The focus of \cite{covert} is
to study the cardinality of the volume set
\[ \tmop{vol} (E) =\{\tmop{vol} (\tmmathbf{x}^1, \ldots, \tmmathbf{x}^d) :
   \tmmathbf{x}^j \in \mathcal{E}\}, \]
where $\mathcal{E}$ is a large subset of $\mathbbm{F}_q^d$. A subset $\mathcal{E} \subset \mathbbm{F}_q^3$ is called a product-like set if  $|\mathcal{E} \cap \mathcal{H}_n | \lesssim
  |\mathcal{E}|^{n / 3}$ for any $n$-dimensional subspace $\mathcal{H}_n \subset \mathbbm{F}_q^3$. Covert et al.
\cite{covert} showed that

\begin{theorem}\label{t3}
  (\cite[Theorem 2.6]{covert}) Suppose that $\mathcal{E} \subseteq \mathbbm{F}_q^3$
  is product-like and $t \in \mathbbm{F}_q^{*}$, then
  \[ |\{\tmop{vol} (\tmmathbf{x}^1, \tmmathbf{x}^2, \tmmathbf{x}^3) = t :
     \tmmathbf{x}^j \in \mathcal{E}\}| = (1 + o (1))
     \frac{|\mathcal{E}|^3}{q}, \]
  if $|\mathcal{E}| = \omega( q^{15 / 8})$.
\end{theorem}

Note that Theorem 2.6 in \cite{covert} only states that $\mathbbm{F}_q^{\ast} \subseteq
\tmop{vol} (\mathcal{E})$ if $|\mathcal{E}| \gg q^{15 / 8}$ but the given
proof in \cite{covert} indeed implies Theorem \ref{t3} above. We will use the geometry incidence machinery developed in that paper \cite{covert} and some properties of non-singular matrices to obtain the following asymptotic result for higher dimensional
cases.

\begin{theorem}\label{m1} For $t \in \mathbbm{F}_q^{\ast}$, $d \geqslant 2$ and  $\mathcal{A} \subset \mathbbm{F}_q$, we have
  \[ N_d (\mathcal{A}; t) = (1 + o (1)) \frac{|\mathcal{A}|^{d^2}}{q}, \]
  if $|\mathcal{A}| = \omega ( q^{\frac{d}{2d-1}} + q^{\frac{d^2-d+4}{2(d^2-d+2)}} )$. 
\end{theorem}

Note that the bound in Theorem \ref{m1} is $|\mathcal{A}| = \omega(q^{\frac{d}{2d-1}})$ if $d \geq 4$ and $|\mathcal{A}| =  \omega (q^{\frac{d^2-d+4}{2(d^2-d+2)}})$ if $d = 2, 3$. When $d=3$, Theorem \ref{m1} matches with the bound in Theorem \ref{t3}, however the later one holds for more general sets. Covert et al. \cite{covert} did not extend their result (Theorem \ref{t3} above) to higher dimensional cases as their focus is the function $|\tmop{vol} (\mathcal{E})|$. They instead showed that $| \tmop{vol} (\mathcal{E}) | =\mathbbm{F}_q$ if $\mathcal{E}=\mathcal{A}
\times \mathcal{A} \times \mathcal{A} \times \mathcal{A}$ whenever $|\mathcal{A}| > \sqrt{q}$. It seems that their proof can be extended to higher dimensional cases.

When $q=p$ is a prime and the set $\mathcal{A}$ is an interval $[- H, H] \subset [-(p-1)/2,(p-1)/2]$, using Theorem \ref{t2}, we
obtain a stronger result for $3\times 3$ matrices.

\begin{theorem} \label{m2}
  For $1 \leq H \leq (p-1)/2$ and $t \in \mathbbm{F}_p^{\ast}$, we have
  \[ N_3 ([- H, H] ; t) = (1 + o (1)) \frac{(2H + 1)^{9}}{p}, \]
  if $H = \omega( p^{\frac{3}{5}})$.
\end{theorem}

Note that the implied constants in the symbols $O, o, \Theta, \Omega, \omega$, and $\ll$ may depend on integer parameter $d$. We recall that the notation $U=O(V)$ and $U \ll V$ are equivalent to the assertion that the inequality $|U| \leq cV$ holds for some constant $c>0$. The notation $U=\Omega(V)$ is equivalent to the assertion that $U \geq c|V|$ holds for some constant $c>0$. We say that $U = o(V)$ if $U = O(V)$ but $U \neq \Omega(V)$ and $U = \omega(V)$ if $U = \Omega(V)$ but $U \neq O(V)$.

\section{Some estimates}

\subsection{Geometric Incidence Estimate}

Let $f$ be a complex-valued function on $\mathbbm{F}_q^d$, we define the $r$-norm of $f$ on $\mathbbm{F}_q^d$ by 
\[ \| f \|_r = \left( \sum_{\tmmathbf{x} \in \mathbbm{F}_q^d} |f(\tmmathbf{x})|^r\right)^{1/r}.\]
The Fourier transform of $f$ on $\mathbbm{F}_q^d$ with respect to a non-trivial principal additive character $\chi$ on $\mathbbm{F}_q$ is given by
\[\hat{f}(\tmmathbf{m}) = q^{-d} \sum_{\tmmathbf{x} \in \mathbbm{F}_q^d} f(\tmmathbf{x})\chi(-\tmmathbf{x} \cdot \tmmathbf{m}).\]
One of our main tools is the following geometric incidence estimate which was
developed and used in \cite{covert} (see also \cite{hart-iosevich,hart} for earlier versions of this estimate).

\begin{theorem} \label{t4}
  (\cite[Theorem 2.1]{covert}) Let $B (\cdot, \cdot)$ be any nondegenerate bilinear
  form in $\mathbbm{F}_q^d$. Let
  \[ \nu (t) = \sum_{B (\tmmathbf{x}, \tmmathbf{y}) = t} f (\tmmathbf{x}) g
     (\tmmathbf{y}), \]
  where $f, g$ are non-negative functions on $\mathbbm{F}_q^d$. Then
  \[ \nu (t) = q^{- 1} \|f\|_1 \|g\|_1 + R (t), \]
  where
  \begin{equation}\label{e1}
    |R (t) | \leqslant q^{\frac{d - 1}{2}} \|f\|_2 \|g\|_2,
  \end{equation}
  if $t \neq 0$. Moreover, if $(0, \ldots, 0) \nin$support$(f) \equiv E$, then
  \begin{equation}\label{e2} \sum_{t \in \mathbbm{F}_q} \nu^2 (t) \leqslant q^{- 1} \|f\|_2^2 \cdot
     |E| \cdot \|g\|_1^2 + q^{2 d - 1} \|f\|_2^2 \sum_{\tmmathbf{k} \neq (0,
     \ldots, 0)} | \hat{g} (\tmmathbf{k}) |^2 |E \cap l_{\tmmathbf{k}} |, \end{equation}
  where
  \[ l_{\tmmathbf{k}} =\{t\tmmathbf{k}: t \in \mathbbm{F}_q^{\ast} \}. \]
\end{theorem}

Note that the proof of Theorem 2.1 in \cite{covert} is given only in the case of dot
product, but it goes through without any essential changes if the dot product
is replaced by any non-degenerate bilinear form.

Theorem \ref{t4} has several applications in additive combinatorics (see \cite{covert,hart-iosevich,hart}). We present here another application of this theorem to the problem of finding three term arithmetic progression in productsets over finite fields. Using multiplicative character sums, Shparlinski \cite{igor3} showed that for any integer $k$ with $p>k\geq3$, where $p$ is the characteristic of $\mathbbm{F}_q$, and any two sets $\mathcal{A}, \mathcal{B} \subset \mathbbm{F}_q$ with
\[|\mathcal{A}| | \mathcal{B}| \geq (k-1)^{2/(k-1)}q^{2-1/(k-1)}, \]
the productset $\mathcal{A}\mathcal{B}$ contains a $k$-term artihmetic progression. He asked if one can relax the condition $k < p$. We give an affirmative answer for this question in the easiest case, $k=3$. It is enough to show that the following equation has solution 
\begin{equation}\label{v} x_0y_0 + x_2y_2 = 2x_1y_1,\,\,\,x_i \in \mathcal{A}, y_i \in \mathcal{B}_i,\end{equation}
has a solution given that $x_0y_0,x_2y_2 \neq x_1y_1$. Fix some $x_1 \in \mathcal{A}, y_1 \in \mathcal{B}$ such that $x_1 y_1 \neq 0$. From (\ref{e1}), the number of quadtuples $(x_0,y_0,x_2,y_2)$ satisfying Eq. (\ref{v}) is at least
\[\frac{|\mathcal{A}|^2 |\mathcal{B}|^2}{q} - \sqrt{q}|\mathcal{A}| |\mathcal{B}|.\] 
Besides, the number of quadtuples $(x_0,y_0,x_2,y_2)$ with $x_0y_0 = x_2y_2 = x_1y_1$ is bounded by $|\mathcal{A}| |\mathcal{B}|$ (as for each $(x_0,y_2) \in \mathcal{A} \times \mathcal{B}$, we have at most one choice for $(y_0,x_2)$). Therefore, the productset $\mathcal{A}\mathcal{B}$ contains a $3$-term artihmetic progression if $|\mathcal{A}| |\mathcal{B}|>q(\sqrt{q}+1)$. Note that for $k=3$, the question of \cite{igor3} is indeed a question about vanishing bilinear forms, so there is no surprise that it admits a different approach using exponential sums, which however is not likely to help for $k>3$. 
\subsection{Recursive estimates}

Let $N_d (\mathcal{A}; t)$ be the number of $d \times d$ matrices with
entries in $\mathcal{A}$ having determinant $t$. The following theorem says that $N_d (\mathcal{A}; t)$ can be bounded by $N_{d-1} (\mathcal{A}; l)$'s.

\begin{theorem} \label{m3}
  For any $t \in \mathbbm{F}_q^{\ast}$ then
  \[ \left| N_d (\mathcal{A}; t) - \frac{|\mathcal{A}|^{d^2}}{q} \right|^2
     \leqslant q^{d - 1} |\mathcal{A}|^{2 d - 1} (1 + o (1)) \sum_{l \in
     \mathbbm{F}_q^{\ast}} N_{d - 1}^2 (\mathcal{A}; l) . \]
\end{theorem}

\begin{proof}
  For any $M \in \mathcal{M}_{d - 1, d} (\mathcal{A})$, let
    $\tmmathbf{m}_i$ be the $i^{\tmop{th}}$ column of $M$ and $M_i$ be the $(d - 1)
    \times (d - 1)$ minor of $M$ by deleting $\tmmathbf{m}_i$. Define
    \[ v (M) = ((- 1)^i \det (M_i))_{1 \leqslant i \leqslant d} \in \mathbbm{F}_q^d. \]
  For any $\tmmathbf{x}= (x_1, \ldots, x_d) \in \mathbbm{F}_q^d$, let $f
  (\tmmathbf{x}) \assign \mathcal{A}^d (\tmmathbf{x}) =\mathcal{A}(x_1) \ldots
  \mathcal{A}(x_d)$ where $\mathcal{A}(\cdot)$ is the characteristic function of the set $\mathcal{A}$, and define
  \[ g (\tmmathbf{x}) \assign |\{M \in \mathcal{M}_{d - 1, d} (\mathcal{A})
     : v (M) =\tmmathbf{x}\}|. \]
  It follows that
  \[ N_d (\mathcal{A}; t) = \sum_{\tmmathbf{x} \cdot \tmmathbf{y}= t} f
     (\tmmathbf{x}) g (\tmmathbf{y}) . \]
  We have $\|f\|_1 =\|f\|_2 = |\mathcal{A}|^d$ and $\|g \|_1 =
  |\mathcal{A}|^{(d - 1) d}$. From (\ref{e1}), we have
  \begin{equation}\label{e8} \left| N_d (\mathcal{A}; t) - \frac{|\mathcal{A}|^{d^2}}{q} \right|^2
     \leqslant q^{d - 1} |\mathcal{A}|^d \|g \|_2^2 . \end{equation}
  Now, we estimate $\|g \|_2^2$. Note that $\tmmathbf{x} \cdot \tmmathbf{y}=
  t \in \mathbbm{F}_q^{*}$ so $\tmmathbf{y} \neq (0, \ldots, 0)$. Therefore
  \begin{eqnarray}
    \|g \|_2^2 & = & \sum_{\tmmathbf{y} \neq (0, \ldots, 0)} g^2
    (\tmmathbf{y}) \nonumber\\
    & = & \sum_{i = 1}^d \sum_{y_i \in \mathbbm{F}_q^{\ast}} \sum_{y_j \in
    \mathbbm{F}_q, j > i} g^2 (0, \ldots, 0, y_i, \ldots, y_d) . \label{e3}
  \end{eqnarray}
  We need the following lemma.
  
  \begin{lemma} \label{l1}
    For any $1 \leqslant i \leqslant d$, then
    \[ \sum_{y_i \in \mathbbm{F}_q^{\ast}} \sum_{y_j \in \mathbbm{F}_q, j > i}
       g^2 (0, \ldots, 0, y_i, \ldots, y_d) \leqslant |\mathcal{A}|^{d - i}
       \sum_{l \in \mathbbm{F}_q^{\ast}} N_{d - 1}^2 (\mathcal{A}; l) . \]
  \end{lemma}
  
  \begin{proof}
    (of the lemma) For any $M \in \mathcal{M}_{d - 1, d} (\mathcal{A})$, let
    $\tmmathbf{m}_i$ be the $i^{\tmop{th}}$ column of $M$ and $M_i$ be the $(d - 1)
    \times (d - 1)$ minor of $M$ by deleting $\tmmathbf{m}_i$.
    For any fixed $y_i,\ldots,y_d \in \mathbbm{F}_q$ and $M_i \in \mathcal{M}_{d - 1, d - 1} (\mathcal{A})$ with
    $\det (M_i) = (- 1)^i y_i \in \mathbbm{F}_q^*$. Let
    \[ \tmmathbf{y}= \frac{1}{\det (M_i)} (0, \ldots, 0, y_{i+1}, \ldots, y_d)^t
       \in \mathbbm{F}_q^d. \]
    We have
    \[ v (M) = ((- 1)^i \det (M_i))_{1 \leqslant i \leqslant d} = (0,\ldots,0,y_i,\ldots,y_d) . \]
    Hence, by Cramer's rule and the non-singularity of $M_i$, we have
    \begin{equation}\label{e4} M_i \tmmathbf{y}=\tmmathbf{m}_i . \end{equation}
    So there is at most one possibility of $\tmmathbf{m}_i$ for each fixed
    $y_i,\ldots,y_d$ and $M_i$. This implies that
    \begin{equation}\label{e5} g (0, \ldots, 0, y_i, \ldots, y_d) \leqslant N_{d - 1} (\mathcal{A};
       (- 1)^i y_i), \end{equation}
       for any $y_i \in \mathbbm{F}_q^*$. Since $\det (M_i) = (-1)^i y_i \in \mathbbm{F}_q^*$, we can write (\ref{e4}) as
    \[ (0, \ldots, 0, y_{i+1}, \ldots, y_d)^t = \det (M_i) M_i^{- 1}
       \tmmathbf{m}_i . \]
    By Gaussian  elimination, we can remove all nonzero entries under the main diagonal in the first $i-1$ rows of $\det(M_i) M_i^{-1}$. Since $\tmmathbf{m}_i \in \mathcal{A}^{d-1}$, for any fixed $M_i$, there are at most $|\mathcal{A}|^{d - i}$ possibilities for $(y_{i
    + 1}, \ldots, y_d)$. This implies that, for any $y_i \in \mathbbm{F}_q^*$ then
    \begin{equation}\label{e6} \sum_{y_j \in \mathbbm{F}_q, j > i} g (0, \ldots, 0, y_i, \ldots,
       y_d) \leqslant |\mathcal{A}|^{d - i} N_{d - 1} (\mathcal{A}; (- 1)^i
       y_i) . \end{equation}
    If $0 \leqslant x, y \leqslant A$, then $x^2 + y^2 \leqslant (\max \{A, x
    + y\})^2 + (x + y - \max \{A, x + y\})^2$. Thus, from (\ref{e5}) and (\ref{e6}), we have
    \[ \sum_{y_j \in \mathbbm{F}_q, j > i} g^2 (0, \ldots, 0, y_i, \ldots,
       y_d) \leqslant |\mathcal{A}|^{d - i} N_{d - 1}^2 (\mathcal{A}; (- 1)^i
       y_i) . \]
    Taking sum over all $y_i \in \mathbbm{F}_q^{\ast}$, the lemma follows.
  \end{proof}
  
  From (\ref{e3}) and Lemma \ref{l1} , we have
  \begin{equation}\label{e7} \|g \|_2^2 \leqslant (|\mathcal{A}|^{d - 1} + \ldots + 1) \sum_{l \in
     \mathbbm{F}_q^{\ast}} N_{d - 1}^2 (\mathcal{A}; l) = |\mathcal{A}|^{d -
     1} (1 + o (1)) \sum_{l \in \mathbbm{F}_q^{\ast}} N_{d - 1}^2
     (\mathcal{A}; l) . \end{equation}
  The theorem follows immediately from (\ref{e8}) and (\ref{e7}). 
\end{proof}

\begin{theorem} \label{m4}
  For any $d \geqslant 2$, then
  \[ \sum_{l \in \mathbbm{F}_q^*} N_d^2 (\mathcal{A}; l) \leqslant
     (1+o(1))\frac{|\mathcal{A}|^{2 d^2}}{q} + q^{d - 1} |\mathcal{A}|^{2 d} (1 + o
     (1)) \sum_{t \in \mathbbm{F}_q^*} N^2_{d - 1} (\mathcal{A}; t) . \]
\end{theorem}

\begin{proof}
  Similarly as in the proof of Theorem \ref{m3}, for any $\tmmathbf{x}= (x_1,
  \ldots, x_d) \in \mathbbm{F}_q^d$, let $f (\tmmathbf{x}) \assign
  \mathcal{A}^d (\tmmathbf{x}) =\mathcal{A}(x_1) \ldots \mathcal{A}(x_d)$, and
  define
  \[ g (\tmmathbf{x}) \assign |\{M \in \mathcal{M}_{d - 1, d} (\mathcal{A})
     : v (M) =\tmmathbf{x}\}|. \]
  Let $f_0(\tmmathbf{x}) = f(\tmmathbf{x})$, $g_0(\tmmathbf{x}) = g(\tmmathbf{x})$ if $\tmmathbf{x} \neq (0,\ldots,0)$ and $f_0(\tmmathbf{x}) = g_0(\tmmathbf{x}) = 0$ otherwise. Then 
  \[ N_d (\mathcal{A}; t) = \sum_{\tmmathbf{x} \cdot \tmmathbf{y}= t} f_0
     (\tmmathbf{x}) g_0 (\tmmathbf{y}), \]
     if $t \in \mathbbm{F}_q^*$.
  Since $(0, \ldots, 0) \nin$support($f_0$)$\equiv \mathcal{E} \subseteq \mathcal{A}^d$, from (\ref{e2}) and Plancherel's theorem, we have
  \begin{eqnarray}
    & & \sum_{l \in \mathbbm{F}_q^*} N_d^2 (\mathcal{A}; l)  \leqslant \sum_{t\in \mathbbm{F}_q} \left(\sum_{\tmmathbf{x} \cdot \tmmathbf{y}= t} f_0
     (\tmmathbf{x}) g_0 (\tmmathbf{y})\right)^2  \nonumber\\
    & \leqslant & q^{- 1} \|f_0\|_2^2 \cdot |\mathcal{E}| \cdot \|g_0 \|_1^2 + q^{2 d - 1}
    \|f_0\|_2^2 \sum_{\tmmathbf{k} \neq (0, \ldots, 0)} | \widehat{g_0}
    (\tmmathbf{k}) |^2 |\mathcal{E} \cap l_{\tmmathbf{k}} |  \nonumber\\
    & \leqslant & \frac{|\mathcal{A}|^{2 d^2}}{q} + q^{2 d - 1}
    |\mathcal{A}|^{d + 1} q^{- d} \sum_{\tmmathbf{y} \in \mathbbm{F}_q^d}
    g_0^2 (\tmmathbf{y})  \nonumber\\
    & = & (1+o(1))\frac{|\mathcal{A}|^{2 d^2}}{q} + q^{d - 1} |\mathcal{A}|^{d + 1}
    \sum_{\tmmathbf{y} \neq (0, \ldots, 0)} g^2 (\tmmathbf{y}), \label{e10}
  \end{eqnarray}
  since $|\mathcal{E} \cap l_{\tmmathbf{k}} | \leqslant |\mathcal{A}|$ for any
  $\tmmathbf{k} \neq (0, \ldots, 0)$. From (\ref{e7}), we have
  \begin{equation} \label{e11} \sum_{\tmmathbf{y} \neq (0, \ldots, 0)} g^2 (\tmmathbf{y}) \leqslant (1
     + o (1)) |\mathcal{A}|^{d - 1} \sum_{l \in \mathbbm{F}_q^{\ast}} N_{d -
     1}^2 (\mathcal{A}; l) . \end{equation}
  The theorem follows from (\ref{e10}) and (\ref{e11}). 
\end{proof}

\section{Distribution of determinant}

\subsection{Arbitrary sets (Proof of Theorem \ref{m1})}

From Theorem \ref{m4}, we have the following corollary.

\begin{corollary}\label{c1}
  For $\mathcal{A} \subseteq \mathbbm{F}_q$ and $d \geqslant 2$, we have
  \[ \sum_{l \in \mathbbm{F}_q^{\ast}} N_d^2 (\mathcal{A}; l) = O \left(
     q^{-1}|\mathcal{A}|^{2 d^2} + q^{\frac{d (d - 1)}{2}}
     |\mathcal{A}|^{d (d + 1) - 1} \right) . \]
\end{corollary}

\begin{proof}
  The proof proceeds by induction. For the base csae $d = 2$, it follows from
  Theorem \ref{m4} that
  \begin{eqnarray}
    \sum_{l \in \mathbbm{F}_q^{\ast}} N_2^2 (\mathcal{A}; l) & \leqslant & (1
    + o (1)) \frac{|\mathcal{A}|^8}{q} + (1 + o (1)) q|\mathcal{A}|^4 \sum_{t \in \mathbbm{F}_q^{*}}
    N_1^2 (\mathcal{A}; t) \nonumber\\
    & = & O \left( q^{-1}|\mathcal{A}|^8 + q|\mathcal{A}|^5 \right) .
    \nonumber
  \end{eqnarray}
  Suppose that the corollary holds for $d - 1$, we show that it also holds for
  $d$. By induction hypothesis, we have
  \begin{equation} \label{e-det1}
    \sum_{t \in \mathbbm{F}_q^{\ast}} N_{d - 1}^2 (\mathcal{A}; t) = O \left(
    q^{-1}|\mathcal{A}|^{2 (d - 1)^2} + q^{\frac{(d - 1) (d - 2)}{2}}
    |\mathcal{A}|^{(d - 1) d - 1} \right) .
  \end{equation}
  Theorem \ref{m4} implies that
  \begin{eqnarray}
    \sum_{l \in \mathbbm{F}_q^{\ast}} N_2^2 (\mathcal{A}; l) & \leqslant & (1
    + o (1)) \frac{|\mathcal{A}|^{2 d^2}}{q} + (1 + o (1)) q^{d - 1}
    |\mathcal{A}|^{2 d} \sum_{l \in \mathbbm{F}_q^{\ast}} N_{d - 1}^2
    (\mathcal{A}; l) \nonumber\\
    & = & O \left( q^{-1}|\mathcal{A}|^{2 d^2} + q^{d - 2}
    |\mathcal{A}|^{2 d^2 - 2 d + 2} + q^{\frac{d (d - 1)}{2}} |\mathcal{A}|^{d
    (d + 1) - 1} \right) \nonumber\\
    & = & O \left( q^{-1}|\mathcal{A}|^{2 d^2} + q^{\frac{d (d - 1)}{2}}
    |\mathcal{A}|^{d (d + 1) - 1} \right), \nonumber
  \end{eqnarray}
  where the second line follows from (\ref{e-det1}) and the last line follows from
  \[ q^{d - 2} |\mathcal{A}|^{2 d^2 - 2 d + 2} = O \left(
     q^{-1}|\mathcal{A}|^{2 d^2} + q^{\frac{d (d - 1)}{2}}
     |\mathcal{A}|^{d (d + 1) - 1} \right) . \]
  This completes the proof of the corollary.
\end{proof}

We are now ready to give a proof of Theorem \ref{m1}. It follows from Theorem \ref{m3}
and Corollary \ref{c1} that
\begin{eqnarray}
  \left| N_d (\mathcal{A}; t) - \frac{|\mathcal{A}|^{d^2}}{q} \right|^2 &
  \leqslant & q^{d - 1} |\mathcal{A}|^{2 d - 1} (1 + o (1)) \sum_{l \in
  \mathbbm{F}_q^{\ast}} N_{d - 1}^2 (\mathcal{A}; l) \nonumber\\
  & = & O \left( q^{d - 2} |\mathcal{A}|^{2 d^2 - 2 d + 1} + q^{\frac{d (d -
  1)}{2}} |\mathcal{A}|^{d (d + 1) - 2} \right) \nonumber\\
  & = & o \left( q^{-2}|\mathcal{A}|^{2 d^2}\right), \nonumber
\end{eqnarray}
given that
\[ |\mathcal{A}| = \omega( q^{\frac{d}{2 d - 1}} + 
   q^{\frac{d^2 - d + 4}{2 (d^2 - d + 2)}} ) . \]
This completes the proof of the theorem.

\subsection{Intervals (Proof of Theorem \ref{m2}) }

It follows from Theorem \ref{t2} that
\begin{eqnarray}
  \sum_{l \in \mathbbm{F}_q^{\ast}} N_2^2 ([- H, H] ; l) & \leqslant & (p - 1)
  \left( \frac{(2 H + 1)^4}{p} + O (H^2 p^{o (1)}) \right)^2 \nonumber\\
  & = & O \left( p^{-1}H^8 + p^{1 + o (1)} H^4 \right) . \label{e-det3}
\end{eqnarray}
From Theorem \ref{m3} and (\ref{e-det3}), we have
\begin{eqnarray}
  \left| N_3 ([- H, H] ; t) - \frac{(2 H + 1)^9}{p} \right|^2 & \leqslant &
  p^2 (2 H + 1)^5 (1 + o (1)) \sum_{l \in \mathbbm{F}_q^{\ast}} N_2^2 ([- H ;
  H) ; l) \nonumber\\
  & = & O \left( p H^{13} + p^{3 + o (1)} H^9 \right) . \nonumber
\end{eqnarray}
This implies that
\[ N_3 ([- H, H] ; t) = (1 + o (1)) \frac{(2 H + 1)^9}{p} \]
if $H = \omega (p^{3 / 5})$, completing the proof of Theorem \ref{m2}.

\section{Remarks}

Note that the quantity
\[ \sum_{l \in F_q^{\ast}} N_d^2 (\mathcal{A}; l) \]
is equal to the number of matrices $M, N$ with entries from $\mathcal{A}$ such
that $M N^{- 1} \in \tmop{SL}_q (d)$. Let $S_d (\mathcal{A})$ denotes this
quantity, it follows from Corollary \ref{c1} that
\[ S_d (\mathcal{A}) = O \left( q^{-1}|\mathcal{A}|^{2 d^2} + q^{\frac{d
   (d - 1)}{2}} |\mathcal{A}|^{d (d + 1) - 1} \right) . \]
Similarly, one can get an estimate of this type for $S_d ([- H
; H])$ from Theorem \ref{t2} and Theorem \ref{m4}.

\begin{corollary}\label{c2}
  Suppose that $p$ is a prime and $1 \leqslant H \leqslant (p - 1) / 2$. We
  have
  \[ S_d([- H, H]) = O \left(
     p^{-1}H^{2 d^2} + p^{\frac{d (d - 1)}{2} + o (1)} H^{d (d + 1) - 2}
     \right) . \]
\end{corollary}

\begin{proof}
  The proof proceeds by induction. For the base case $d = 2$, it follows from
  Theorem \ref{t2} that
  \begin{eqnarray}
    S_2([- H, H]) & \leqslant & (p -
    1) \left( \frac{(2 H + 1)^4}{p} + O (H^2 p^{o (1)}) \right)^2 \nonumber\\
    & = & O \left( p^{-1}H^8 + p^{1 + o (1)} H^4 \right) . \nonumber
  \end{eqnarray}
  Suppose that the corollary holds for $d - 1$, we show that it holds for $d$.
  Theorem \ref{m4} implies that
  \begin{eqnarray}
    S_d([- H, H]) & \leqslant & (1 +
    o (1)) \frac{(2 H + 1)^{2 d^2}}{p} + (1 + o (1)) p^{d - 1} (2 H + 1)^{2 d}
    S_{d-1}([- H, H]) \nonumber\\
    & = & O \left( p^{-1}H^{2 d^2} + p^{d - 2} H^{2 d^2 - 2 d + 2} +
    p^{\frac{d (d - 1)}{2} + o (1)} H^{d (d + 1) - 2} \right) \nonumber\\
    & = & O \left( p^{-1}H^{2 d^2} + p^{\frac{d (d - 1)}{2} + o (1)} H^{d
    (d + 1) - 2} \right), \nonumber
  \end{eqnarray}
  where the second line follows from the induction hypothesis and the last
  line follows from
  \[ p^{d - 2} H^{2 d^2 - 2 d + 2} = O \left( p^{-1}H^{2 d^2} + p^{\frac{d
     (d - 1)}{2} + o (1)} H^{d (d + 1) - 2} \right) . \]
  This completes the proof of the corollary.
\end{proof}

It has been pointed out by the referee that the bound in Corollary \ref{c2} can be
improved for small value of $H$.

\begin{lemma}\label{r-det}
  Suppose that $q = p$ is a prime and $H = O (p^{1 / 2})$, then
  \[ S_d ([- H, H]) = O \left( p^{\frac{d (d - 1)}{2} - 1} H^{d (d + 1) + o
     (1)} \right) . \]
\end{lemma}

\begin{proof}
  The proof proceeds by induction. For the base case $d = 2$ we have
  \[ S_2 ([- H, H]) \leqslant |\{x_1 y_2 - x_2 y_1 \equiv u_1 v_2 - u_2 v_1 
     \: (\tmop{mod} p)\}| \]
  where all variables are in $[- H, H]$. For each choice of $x_2, y_1, u_1,
  v_1, u_2$ and $v_2$, we get $x_1 y_2 \equiv a$ (mod $p$) for some $a \in
  \mathbbm{F}_q$.
  
  If $a = 0$, then there are $O (H^5)$ posibilities for $x_2, y_1, u_1, v_1,
  u_2, v_2$ and $O (H)$ posibilities for $x_1, y_2$. If $a \in
  \mathbbm{F}_q^{\ast}$, then the arithmetic progression $z \equiv a$ (mod
  $p$) contains $O (H^2 / p + 1)$ elements $|z| \leqslant H^2$. Since $z \neq
  0$ for all of them, $z = x_1 y_2$ has $H^{o (1)}$ solutions. Putting
  everything together, we get
  \[ S_2 ([- H, H]) = O (H^6 + H^6 (H^2 / p + 1) H^{o (1)}) = O (H^{6 + o
     (1)}) \]
  if $H = O (p^{1 / 2})$. Suppose that the lemma holds for $d - 1$, we show
  that it holds for $d$. Theorem \ref{m4} implies that
  \begin{eqnarray}
    S_d ([- H, H]) & = & O (H^{2 d^2} / p + p^{d - 1} H^{2 d} S_{d - 1} ([- H
    ; H])) \nonumber\\
    & = & O (H^{2 d^2} / p + p^{\frac{d (d - 1)}{2} - 1} H^{d (d + 1) + o
    (1)}) \nonumber\\
    & = & O (p^{\frac{d (d - 1)}{2} - 1} H^{d (d + 1) + o (1)}), \nonumber
  \end{eqnarray}
  where the second line follows from the induction hypothesis and the last
  line follows from $H = O (p^{1 / 2})$. This completes the proof of the
  lemma.
\end{proof}

\section*{Acknowledgement}
The author would like to thank an anonymous referee for detailed comments and suggestions, especially for pointing out Lemma \ref{r-det}. He also wants to thank Dang Phuong Dung for carefully reading the manuscript.

\end{document}